\documentclass[11pt,a4paper]{amsart}
\usepackage{latexsym}
\usepackage{amssymb}
\usepackage{calc}
\usepackage{layout}
\usepackage[bookmarks=true]{hyperref}   
\usepackage{xspace}
\usepackage[all]{xy}        
\CompileMatrices            
\UseTips                    

\theoremstyle{plain}
\newtheorem{theorem}{Theorem}[section]
\newtheorem*{theorem*}{Theorem}
\newtheorem{proposition}[theorem]{Proposition}
\newtheorem{corollary}[theorem]{Corollary}
\newtheorem{lemma}[theorem]{Lemma}

\theoremstyle{definition}

\theoremstyle{remark}
\newtheorem{remark}[theorem]{Remark}

\numberwithin{equation}{section}

\newcommand{\F}[1][]{F^{#1*}}                  
\newcommand{\D}[1][]{D^{(#1)}}                 

\newcommand{\enm}[1]{\ensuremath{#1}}          
\newcommand{\op}[1]{\operatorname{#1}}

\newcommand{\ZZ}{\enm{\mathbb{Z}}}

\renewcommand{\phi}{\varphi}        
\renewcommand{\theta}{\vartheta}
\renewcommand{\epsilon}{\varepsilon}

\newcommand{\Ann}{\op{Ann}}         

\newcommand{\Spec}{\op{Spec}}

\newcommand{\Hom}{\op{Hom}}

\newcommand{\End}{\op{End}}

\newcommand{\id}{\op{id}}

\newcommand{\tensor}{\otimes}         

\renewcommand{\to}[1][]{\xrightarrow{\ #1\ }}

\newcommand{\usc}[1][m]{\underline{\phantom{#1}}}

\newcommand{\defeq}{\stackrel{\scriptscriptstyle \op{def}}{=}}

\newcommand{\cf}{\textit{cf.}\ }

\begin{document}
\title[]{$D$--module generation in positive characteristic via Frobenius descent}
\author{Manuel Blickle}
\address{Universit\"at Essen, FB6 Mathematik, 45117 Essen,
Germany} \email{manuel.blickle@uni-essen.de}
\urladdr{\url{www.mabli.org}}


\maketitle

\begin{abstract}
    In this note I show that for most $F$--finite regular rings $R$ of positive
    characteristic the localization $R_f$ at an element $f\in R$ is generated
    by $f^{-1}$ as a $D_R$--module. This generalizes and gives an alternative
    proof of te results in \cite{MontLyub} where the result is proven for the polynomial
    ring, thereby answering questions raised in \cite{MontLyub} affirmatively.
    The proof given here is a surprisingly simple application of
    Frobenius descent, a brief but thorough discussion of which is also
    included. Furthermore I show how essentially the same technique yields a
    quite general criterion for obtaining $D_R$--module generators of a unit
    $R[F]$--module.
\end{abstract}

\section{$R_f$ is $D_R$--generated by $f^{-1}$}

Throughout this paper $R$ will denote a noetherian regular ring. This note was
created after hearing about the result of Alvarez~Montaner and Lyubeznik in
\cite{MontLyub}. For the case of $R=k[x_1,\ldots,x_n]$ they show the following
theorem:
\begin{theorem}\label{thm.Main1}
    Let $R$ be a regular $F$--finite ring of positive characteristic,
    which is essentially of finite type over a regular local ring.
    Let $f \in R$ be a nonzero element. Then the $D_R$--module
    $R_f$ is generated by $f^{-1}$.
\end{theorem}
What is surprising about this result is that in characteristic zero it is
incorrect. There, the $D_R$--generation of $R_f$ is governed by the
Bernstein-Sato polynomial of $f$, and therefore reflects the geometry of the
hypersurface defined by $f=0$. Thus this is another instance where
$D_R$--modules appear coarser in positive characteristic than in characteristic
zero, \cf for example \cite{Bli.int,LauRaPedFuch-GlobalFreg}.

The proof of Theorem \ref{thm.Main1} given here relies on two results which are
central to the study of differential operators in positive characteristic. The
first one is the fact that the localization $R_f$ has finite length as a
$D_R$--module, which is implied by the following more general result of
Lyubeznik.

\begin{theorem}[\protect{\cite[Theorem 5.7]{Lyub}}]
Let $R$ be regular, $F$--finite and essentially of finite type over a
$F$--finite ring. Let $M$ be a finitely generated unit $R[F]$--module. Then $M$
has finite length as a $D_R$--module.
\end{theorem}

A \emph{unit $R[F]$--module}\footnote{What we call here a \emph{finitely
generated unit $R[F]$--module} is called an \emph{$F$--finite $R$--module} in
\cite{Lyub}} is an $R$--module $M$ together with an isomorphism $\theta_M: F^*M
\to M$. Note that $F: \Spec R \to \Spec R$ is the absolute Frobenius map, which
is the identity on the underlying topological space and the $p$th power map on
the structure sheaf. In particular on global sections this is just the
Frobenius map $F: R \to R$ raising $r$ to $r^p$, which, by abuse, is also
denoted by the letter $F$. The assumption that $R$ be $F$--finite just means
that the Frobenius is a finite map.

The functors $F^*$ and $F_*$ are just pullback and pushforward along the
Frobenius. Concretely, one can think of $F^*M = R {{}_F\negthinspace \tensor}
M$ where the tensor product on the left is via the Frobenius $F:R \to R$.
Similarly $F_*M$ is just $M$ as an abelian group with $R$ structure twisted by
the Frobenius.

By adjointness, a $R$--linear map $\theta: F^*M \to M$ is equivalent to an
$R$--linear map $F_M: M \to F_*M$ which in turn is nothing but a structure of a
module over the noncommutative ring
\[
    R[F] \defeq \frac{R\langle F \rangle}{\{\, r^pF-Fr\,|\, r \in R\, \}}
\]
on $M$. On calls $(M,\theta)$ \emph{finitely generated} if it is finitely
generated as a module over the ring $R[F]$. For $f \in R$ the localization
$R_f$ is a finitely generated unit $R[F]$--module.\footnote{With the
identification $F^*R_f=R {{}_F\negthinspace \tensor} R_f$ the map $R_f \to
F^*R_f$ sending $\frac{a}{b}$ to $ab^{p-1} \tensor \frac{1}{b}$ is inverse to
the natural map $\theta: F^*R_f \to R_f$. The corresponding Frobenius action
$F$ is just raising to the $p$th power. Therefore $R_f$ is generated as an
$R[F]$--module by $f^{-1}$. Thus $R_f$ is a finitely generated unit
$R[F]$--module.}

It was shown in \cite{Lyub} that a (finitely generated) unit $R[F]$--module
$(M,\theta)$ carries a natural structure of a $D_R$--module. Furthermore, the
structural map $\theta: \F M \to M$ is then $D_R$--linear, where the
$D_R$--structure on $F^*M$ is due to the Theorem \ref{thm.FrbDesc} below.

The second crucial ingredient is the so called \emph{Frobenius descent}. For
the convenience of the reader I will include a (very) brief treatment of this
powerful and widely applicable technique at the end of this paper. Most
relevant here is the following consequence:
\begin{theorem}\label{thm.FrbDesc}
Let $R$ be a regular and $F$--finite ring of positive
characteristic. Then the Frobenius functor is an autoequivalence
of the category of $D_R$--modules. In particular, for a
$D_R$--module $M$ the module $F^*M$ carries a natural
$D_R$--module structure.
\end{theorem}

\begin{proof}[Proof of Theorem \ref{thm.Main1}]
For any $D_R$--submodule $M \subseteq R_f$ one identifies the $D_R$--module
$F^*M$ with its isomorphic image in $R_f$ via the natural $D_R$--module
isomorphism $\theta: F^*R_f \to R_f$. Then $F^*M$ is the $D_R$--submodule of
$R_f$ consisting of the elements $rm^p$ for $r \in R$ and $m \in M \subseteq
R_f$.

Let $M=D_Rf^{-1}$ and let me point out that $M \subseteq F^*M$:
Because $F^*M$ is a $D_R$--submodule of $R_f$ (by Frobenius
descent) it is enough to show that $f^{-1} \in F^*M$. Since $F^*M$
(as a submodule of $R_f$) consists precisely of the elements
$rm^p$ for $r \in R$ and $m \in M$ we may use $r=f^{p-1}$ and
$m=f^{-1}$ to conclude that $f^{-1}=rm \in F^*M$.

Now, by repeated application of the Frobenius we get an increasing
chain of $D_R$--submodules of $R_f$:
\begin{equation}\tag{$**$}\label{eq}
    M \subseteq F^*M \subseteq F^{2*}M \subseteq F^{3*}M \subseteq \ldots
\end{equation}

Since $f^{-1} \in M$ it follows that $f^{-p^e}=F^e(f^{-1})$ is an
element of $F^{*e}M$ which shows that the union of the chain must
be all of $R_f$.

Thus it is enough to show that $M=F^*M$ since then the limit
system is constant and $M=R_f$ as claimed. Let us suppose
otherwise, that is assume that the inclusion $M \subsetneq F^*M$
is strict. By Frobenius descent, \emph{all} the inclusions of
(\ref{eq}) must be strict. But this contradicts the finite length
of $R_f$ as a $D_R$--module.
\end{proof}
\begin{remark}
With this result I am able to answer the two (related) questions raised at the
end of \cite{MontLyub}. Firstly, they asked whether their result for the
polynomial ring would also hold for the power series ring; this case is covered
by the above theorem.

Secondly, in their proof is only one step (\cite[Lemma 3.5]{MontLyub}) which
does not hold in the complete case. In fact one easily sees that the truth of
their Lemma 3.5 for a power series ring is in fact equivalent to Theorem
\ref{thm.Main1} for a power series ring. Since the latter was established
above, Lemma 3.5 in \cite{MontLyub} is therefore also valid in the complete
case.
\end{remark}

\section{$D_R$--generators of unit $R[F]$--modules}

The above result also follows from a more general observation.
\begin{theorem}\label{thm.Main2}
Let $R$ be a regular $F$--finite ring of positive characteristic, which is
essentially of finite type over a regular local ring. Let $N$ be a finitely
generated unit $R[F]$--module. Suppose $M \subseteq N$ is a $D_R$--submodule
such that $M \subseteq F^*M$. Then $M$ is a unit $R[F]$--submodule.
\end{theorem}
\begin{proof}
    Once more I identify $F^*M \subseteq F^*N$ with its isomorphic image
    in $N$ via the structural isomorphism $\theta: F^*N \to N$ of the unit
    $R[F]$--module $N$. Then, $M$ being a unit $R[F]$--submodule just means
    that the inclusion $M \subseteq F^*M$ is in fact an equality.
    Assuming otherwise we apply Frobenius
    to the strict inclusion $M \subsetneq F^*M$. By Frobenius descent we
    conclude that all the inclusions $F^{e*}M
    \subsetneq F^{(e+1)*}M$ are strict as well. The resulting strictly increasing infinite
    chain
    \[
        M \subsetneq F^*M \subsetneq F^{2*}M \subsetneq F^{3*}M
        \subsetneq \cdots
    \]
    contradicts the finite length of $N$ as a $D_R$--module.
\end{proof}
This result was inspired by a result in \cite{EmKis.Fcrys}, Proposition 15.3.4,
which (in the notation of Theorem \ref{thm.Main2}) states that if $F^*M
\subseteq M$ then $M$ is also a unit $R[F]$--submodule.

To obtain Theorem \ref{thm.Main1} from this just note that $M =
D_Rf^{-1}$ satisfies $M \subseteq F^*M$ and contains the
$R[F]$--module generator $f^{-1}$ of $R_f$.

\begin{corollary}
With the same assumptions as in Theorem \ref{thm.Main2}, if $n_1,\ldots,n_t$
are generators of a \emph{root}\footnote{An $R$--submodule $N_0$ of a unit
$R[F]$--module $N$ is called a \emph{root}, if $N_0$ is finitely generated as
an $R$--module, $N_0 \subseteq F^*N_0$ and $\bigcup_e^\infty F^{e*}N_0 = N$.
The existence of a root is equivalent to $N$ being finitely generated as a unit
$R[F]$--module.} of the finitely generated unit $R[F]$--module $N$, then
$n_1,\ldots,n_t$ generate $N$ as a $D_R$--module.
\end{corollary}
\begin{proof}
    By Theorem \ref{thm.Main2} it is enough to check that the $D_R$--submodule $M \defeq
    D_R\langle n_1,\ldots,n_t \rangle$ satisfies $M \subseteq
    F^*M$ and contains the $R[F]$--module generators
    $n_1,\ldots,n_t$ of $N$. The second statement is trivial and
    for the first one observes that, by definition of
    root, one can write $n_i = \sum r_j F_N(n_j)$ for some $r_j
    \in R$. Noting that $F_N(n_j) \in F^*M$ we conclude $n_i \in
    F^*M$ for all $i$ as required.
\end{proof}

\section{Frobenius Descent}
Frobenius descent (in the basic form used here) is based on the simple fact
that a ring $R$ is Morita equivalent to the matrix algebra of $n \times n$
matrices with entries in $R$. That is $R$ and $\op{Mat}_{n\times n}(R)$ have
equivalent module categories.

The application of this basic observation to the study of $D_R$--modules in
positive characteristic is very successful as the works of S.P. Smith
\cite{SmithSP.diffop,SmithSP:DonLine}, B. Haastert
\cite{Haastert.DiffOp,Haastert.DirIm} and R. B{\o}gvad \cite{Bog.DmodBorel}
demonstrate.\footnote{A predecessor of it is the so-called Cartier descent as
described, for example, in Katz \cite[Theorem 5.1]{Katz}. It states that $F^*$
is an equivalence between the category of $R$--modules and the category of
modules with integrable connection and $p$--curvature zero. The inverse functor
of $\F$ on a module with connection $(M,\nabla)$ is in this case given by
taking the horizontal sections $\ker \nabla$ of $M$. As an $R$--module with
integrable connection and $p$--curvature zero is nothing but a
$\D[1]_R$--module, Cartier descent is just the case $e=1$ of Proposition
\ref{prop.FrobDesc}.} The ultimate generalization is the treatment of Berthelot
\cite{Ber.FrobDesc}.

The way this Morita equivalence enters the picture comes from the
fact that in characteristic $p > 0$, the ring of differential
operators of an $F$--finite ring $R$ is the union
\[
    D_R = \bigcup_e \D[e]_R
\]
where $\D[e]_R = \End_{R^{p^e}}(R)$ consist of endomorphisms $\phi \in
\End_\ZZ(R)$ which are linear over the subring $R^{p^e}$ of $p$th powers of $R$
(see for example \cite{Yeku} or \cite[Chapter 3.1]{Bli.PhD}). Replacing the
$R^{p^e}$ linear inclusion $R^{p^e} \subseteq R$ with the $R$--linear $F: R \to
F^e_*R$ we identify $\D[e]_R$ with $\End_R(F^e_*R)$. If in addition $R$ is
regular, a basic result of Kunz \cite{Kunz} implies that $F^e_*R$ is a locally
free $R$--module of finite rank, thus locally, $\D[e]_R$ is indeed just a
matrix algebra over $R$.

The aim of this section\footnote{This section is an excerpt of Chapter 3.2 of
\cite{Bli.PhD}. We state and proof the basic result but for all the
straightforward (but tedious) comptibilities one has to check we refer to
\cite{Bli.PhD}. Obviously, in this section no originality beyond the exposition
is claimed.} is to make the resulting Morita equivalence between $R$ and
$\D[e]_R=\End_R(F^e_*R)$ explicit. Concretely I want to show that it is induced
by the Frobenius. This is the content of the following proposition.

\begin{proposition}[Frobenius Descent]\label{prop.FrobDesc}
    Let\/ $R$ be regular and\/ $F$--finite. Then\/ $\F[e]$ is an equivalence of
    categories between the category of\/ $R$--modules and the
    category of $\D[e]_R$--modules. Its inverse functor is given by
    \[
        T^{e}(\usc) \defeq \Hom_R(F^e_*R,R) \tensor_{\D[e]_R}
        \usc.
    \]
\end{proposition}
\begin{proof}
Let us view $\Hom_R(F^e_*R,R)$ as an $R$--$\D[e]_R$--bimodule. The action is by
post-- and pre--composition respectively. That is for $\delta \in \D[e]$, $\phi
\in \Hom_R(F^e_*R,R)$ and $r \in R$ the product $r \cdot \phi \cdot \delta$ is
given by the composition
\[
    F^e_*R \to[\delta] F^e_*R \to[\phi] R \to[r\cdot] R.
\]
Also $F^e_*R$ itself is viewed as a $\D[e]_R$--$R$--bimodule. The left
$\D[e]_R$--structure is clear by definition of $\D[e]_R=\End_R(F^e_*R)$ and the
the right $R$--structure is via Frobenius $F: R \to F^e_*R$. Thus, the
Frobenius functor $F^{e*}(\usc)$ can be identified with $F^e_*R \tensor_R \usc$
and thus $F^*M$ naturally carries the structure of a $\D[e]_R$--module. In this
manner it is clear that the described associations are functors between the
claimed categories as they are just tensoring with an appropriate bimodule.

It remains to show that they are canonically inverse to each other. For this
observe that the natural map
\[
    \Phi:\ F^e_*R \tensor_R \Hom_R(F^e_*R,R) \to \D[e]_R
\]
given by sending $a \tensor \phi$ to the composition
\[
    F_*R \to[\phi] R \to[F^e] F_*R \to[a\cdot] F_*R
\]
is an isomorphism (of bi--modules) by the fact that $F^e_*R$ is a locally free
and finitely generated $R$--module ($\Hom\ $ commutes with finite direct sums
in the second argument). Thus $\Phi$ is a natural transformation of $F^e_*R
\tensor_R \Hom_R(F^e_*R,R) \tensor_{\D[e]_R} \usc$ to the identity functor on
$\D[e]_R$--mod.

Conversely it is equally easy to see that the map
\[
    \Psi:\ \Hom_R(F^e_*R,R) \tensor_{\D[e]_R} F^e_*R \to R
\]
given by sending $\phi \tensor a$ to $\phi(a)$ is also an
isomorphism: After a local splitting $\pi$ of $F^e:R \to F^e_*R$
is chosen (it exists by local freeness of $F^e_*R$ over $R$), its
inverse is given by $a \mapsto \pi \tensor F^e(a)$.
\end{proof}
\begin{remark}
It is possible to give a more explicit description of $T^e$ as follows. Let
$J_e$ be the left ideal of $\D[e]_R$ consisting of all operators $\delta$ such
that $\delta(1)=0$. Then
\[
    T^e(M) \cong \Ann_M J_e.
\]
If one has a splitting $\pi_e$ of the Frobenius $F^e:R \to F^e_*R$ there is yet
another description of $T^e$. Note that $(F^e \circ \pi_e)$ is a map
\[
    F^e_*R \to[\pi_e] R \to[F^e] F^e_*R
\]
such that it can be viewed as a differential operator in $\D[e]_R$ and
therefore acts on any $\D[e]_R$--module $M$. One can show that
\[
    T^e(M) \cong (F^e \circ \pi_e)(M)
\]
and that $(F^e \circ \pi_e)(M) \subseteq M$ is independent of the chosen
splitting $\pi_e$. These statements are also verified in \cite[Chapter
2.3]{Bli.PhD}.
\end{remark}

Proposition \ref{prop.FrobDesc} implies that the categories of
$\D[e]_R$--modules for all $e$ are equivalent since each single
one of them is equivalent to $R$--mod. The functor giving the
equivalence between $\D[f]$--mod and $\D[f+e]$--mod is, of course,
$\F[e]$. Concretely, to understand the $\D[f+e]$--module structure
on $\F[e]M$ for some $\D[f]$--module $M$, we write $M \cong
\F[f]N$ for $N=T^f(M)$. Then $\F[e]M=\F[(f+e)]N=R^{(f+e)}\tensor
N$ carries obviously a $\D[f+e]$--module structure with $\delta
\in \D[f+e]$ acting via $\delta \tensor \id_N$.

Since the union $\bigcup \D[e]_R$ is just the ring of differential
operators $D_R$ of $R$ this shows (after the obvious
compatibilities are checked, which is straightforward and carried
out in \cite[Chapter 3.2]{Bli.PhD}) that $\F[e]$ is in fact an
auto--equivalence of the category of $D_R$--modules:
\begin{proposition}\label{prop.FrobDesc4Dmod}
    Let\/ $R$ be regular and\/ $F$--finite. Then\/ $\F[e]$ is an equivalence of
    the category of\/ $D_R$--modules with itself.
\end{proposition}
Finally we specialize to the case that $M$ is a unit $R[F]$--module. Therefore
it naturally carries a $D_R$--module structure \cite{Lyub}, and so does $F^*M$,
via Frobenius descent. The following lemma shows that these structures are
compatible:
\begin{lemma}\label{lem.ThetaDlinear}
    Let\/ $R$ be regular and\/ $F$--finite and\/ $(M,\theta)$ be a unit\/
    $R[F]$--module. Then\/ $\theta: F^*M \to M$ is a map of\/ $D_R$--modules.
\end{lemma}
\begin{proof}
Again I omit the strightforward verification of this and instead refer to
\cite[Chapter 3.2]{Bli.PhD}
\end{proof}

\providecommand{\bysame}{\leavevmode\hbox to3em{\hrulefill}\thinspace}

\end{document}